\documentclass[11pt]{article}

\usepackage[margin=1in]{geometry}
\usepackage{enumerate, amsmath, amssymb,latexsym}
\usepackage{amsthm, color, graphicx, hyperref, amscd}
\usepackage{todonotes, comment}
\usepackage{bbm}
\newcommand\nc{\newcommand}
\nc\DMO{\DeclareMathOperator}
\nc\ignore[1]{}     

\usepackage[normalem]{ulem}

\newcommand{\cP}{\mathcal{P}}

\newcommand{\cL}{\mathcal{L}}

\newcommand{\B}{\mathcal{B}}
\newcommand{\cB}{\mathcal{B}}

\newcommand{\F}{\mathbb{F}}

\newcommand{\PG}{\mathrm{PG}}

\theoremstyle{definition}

\DMO{\Tr}{Tr}        \DMO{\N}{N}    \DMO{\ind}{ind}
\DMO{\AG}{{AG}}                        \DMO{\GL}{{GL}}
\DMO{\PGL}{{PGL}}                     \DMO{\PGaL}{P$\Gamma$L}
\DMO{\Gal}{Gal}                          \DMO{\GaL}{$\Gamma$L}

\nc{\uE}{{\textup E}}                              \nc{\uL}{{\textup L}}
\nc{\uS}{{\textup S}}                              \nc{\Zp}{{\mathbb Z}_p}
\nc{\cR}{{\mathcal R}}                            \nc{\cA}{{\mathcal A}}

\nc{\cN}{{\mathcal N}}                           \nc{\cM}{{\mathcal M}}
\nc{\cG}{{\mathcal G}}                           \nc{\ts}{{\theta_s}}
\nc{\cE}{{\mathcal E}}                            \nc{\cF}{{\mathcal F}}

\nc{\wa}{{\widehat{\alpha}}}                 \nc{\wg}{{\widehat{\gamma}}}
\nc{\cJ}{{\mathcal J}}                            \nc{\cQ}{{\mathcal Q}}
\nc{\Fqt}{\F_{q^t}}                               \nc{\K}{{\mathbb K}}
\nc{\GF}{\hbox{{\rm GF}}}                    \nc{\st}{{\cal S}_t}
\nc{\so}{{\cal S}}                                  \nc{\sk}{{\cal S}_2}
\nc{\ad}{E}                                            \nc{\lmb}{\lambda}
\nc{\la}{\langle}                                     \nc{\ra}{\rangle}

\nc{\ZG}{\mathcal{Z}(\mathrm{\Gamma L})}
\nc{\G}{\mathrm{\Gamma L}}   
\nc{\gen}{\mathrm{gen}}
\nc{\ep}{\epsilon}
\nc{\vv}{\mathbf v}  
\nc{\1}{\mathbbm 1}  
\nc{\hs}{{\hat{\sigma}}}

\newtheorem{theorem}{Theorem}[section]
\newtheorem{problem}[theorem]{Problem}
\newtheorem{lemma}[theorem]{Lemma}
\newtheorem{corollary}[theorem]{Corollary}
\newtheorem{construction}[theorem]{Construction}
\newtheorem{definition}[theorem]{Definition}
\newtheorem{proposition}[theorem]{Proposition}
\newtheorem{result}[theorem]{Result}

\newtheorem{defi}[theorem]{Definition}

\newtheorem{remark}[theorem]{Remark}
\usepackage{xcolor}
\newcommand{\zs}[1]{\textcolor{magenta}{\textbf{#1}}}
\newcommand{\zsst}[1]{%
	\bgroup
	\markoverwith{\textcolor{magenta}{\rule[0.5ex]{2pt}{1.2pt}}}%
	\ULon{#1}%
	\egroup
}

\title{Extremal secant size distributions in projective planes}
\author{Zolt\'an L\'or\'ant Nagy\thanks{
		E\"otv\"os Lor\'and University, Budapest, Hungary. The author is supported by the J\'anos Bolyai Scholarship of the Hungarian Academy of Sciences and by the NRDI EXCELLENCE-24 grant no. 151504 Combinatorics and Geometry. 
		E-mail: {\tt zoltan.lorant.nagy@ttk.elte.hu}}
	\and Zsuzsa Weiner \thanks{Prezi.com, H-1065 Budapest, Nagymező utca 54-56, Hungary;  E-mail: {\tt zsuzsa.weiner@gmail.com}
}}
\date{}
\date{}

\begin{document}
	
	\maketitle
	
	\begin{abstract}
		Given a point set \(S\) in a projective plane of order \(q\), each line $\ell$ determines a secant of size
		\(|S\cap \ell|\). We study how concentrated this secant-size distribution can be.
		We develop a general framework for estimating and, in several cases, determining the extremal number. We point out how the extremal constructions vary depending on the remainder of the size of the set modulo $q$, and classify several extremal families of point sets, building from concurrent lines, Baer {subplanes} and lines in general position. On the lower bound side, we rely on stability results in the Desarguesian case and applications of the standard equalities in general. 
		Finally, we  connect the problem to a novel concept, called  general position on level 2 for point and line sets, which then leads to a Sidon-type problem over finite fields called projective Sidon sets, providing an additional arithmetic model for
		highly concentrated secant-size distributions.
	\end{abstract}
	
	\section{Introduction}
	
	A recurring theme in combinatorics and incidence geometry is to understand how a set is distributed
	with respect to a family of test objects, such as lines, subspaces or, more generally, hyperedges of a hypergraph defined on the ground set. Given an incidence structure---for instance a block design,
	a finite geometry, or a family of lines or hyperplanes---and a subset \(S\) of the point set, one may
	study the multiset of intersection sizes
	\[
	|S\cap B|
	\]
	as \(B\) ranges over the blocks of the structure. A basic question is how concentrated this
	distribution can be: can one intersection size occur on almost all blocks, or must several different
	intersection sizes appear with comparable frequency?
	
	Questions of this type arise naturally in many settings. In Euclidean incidence geometry, or after projective completion in the real projective plane, a prime example is the orchard problem  due to Jackson \cite{Jack}, which  studies  the maximum number of $3$-rich lines determined by an \(n\)-point set. For the state-of-the-art result on this problem, see the beautiful paper of Green and Tao \cite{Green-Tao}. A related problem concerns $k$-rich lines: lines with prescribed multiplicity $k$.
	In fact, a well-known consequence of the famous Szemerédi--Trotter theorem claims that the number of lines with at least $k\ge2$ points in a planar $n$-set of points is at most $O(n^2/k^3+n/k)$ \cite{Pach, Szem-T}.  This provides indirect control on large secant classes via incidence limits or collinearity avoidance, and forces the number of short secants to be large \textit{in every non-collinear} point set, with respect to the number of all $k$-secants ($k>1$). We mention that this phenomenon is far from true in the finite projective setting, see the recent paper of the authors \cite{WZsNZL}. 
	

	
	In this paper we investigate this question in finite projective planes. 
	Let $\Pi_q$ be a projective plane of order $q$, with point set $\cP$ and line set $\cL$, and let
	$N:=q^2+q+1$
	denote the number of points, and also the number of lines, of $\Pi_q$. If the plane is Desarguesian, we denote it by $\PG(2,q)$.
	
	\begin{defi}[$k$-secants] \label{def:k-secants}
		For a point set $S\subseteq \cP$ and a  line $\ell$, we define $n_{\ell}:=|S\cap \ell|$ to be the intersection size of $\ell$ with $S$.
		A line is called a\textit{ $k$-secant of a set} $S$ if $n_\ell=k$. In particular, for $k=0$ and $k=1$, these lines are usually called \textit{skew} lines and\textit{ tangent }lines of $S$, respectively.\\  For an integer $k$, let $$L_k(S):=\{\ell\in \cL: n_{\ell}=k\}$$ denote the set of $k$-secant lines.
	\end{defi}
	
	\begin{defi}
		\label{def: s_m}
		
		Let $s_m(S)$ denote the maximum number of occurrences of any secant size, i.e., $$s_m(S):=\max_k \{|L_k(S)|\}.$$
	\end{defi}
	
	Rather than prescribing a particular secant size, we ask a more global extremal question
	concerning the distribution of the secant sizes.
	
	\begin{problem} \label{maxi}
		Given a positive integer $n$, determine $$\max \{s_m(S): {|S|=n, S\subseteq \cP}\}.$$ Equivalently, determine how concentrated the secant-size distribution of an $n$-set of points can be.
	\end{problem}
	
	It is convenient to measure the number  of atypical lines instead.
	Accordingly, we define
	\[
	M(n,\Pi_q):=N-\max\{s_m(S):\ S\subseteq \cP,\ |S|=n\}.
	\]
	Hence $M(n,\Pi_q)$ is the minimum possible number of non-most-common secants
	among all point sets of size $n$ in a fixed projective plane $\Pi_q$. 
	For the Desarguesian plane, we write $M^D(n,q):=M(n,\PG(2,q))$.

	Our main observation is that the extremal behaviour depends primarily on the residue of $n$ modulo the order of the {plane}, $q$. Writing
	\[
	n=qr+t,\qquad -q/2 < t \le q/2,
	\]
	
	separates the cardinality into a multiple of \(q\) and a centered remainder. The remainder \(t\) turns out to be a useful parameter for the extremal problem. When \(|t|\) is small, natural candidate configurations arise from unions of concurrent lines and small perturbations thereof. When \(q\) is a square, additional families based on Baer subplanes become competitive for certain remainders on the \(\sqrt q\)-scale. Thus the behaviour of \(M(n,\Pi_q)\) depends strongly on the arithmetic position of \(n\) relative to multiples of $q$.
	
	
	The extremal constructions appearing throughout the paper arise from three geometric mechanisms.
	
	First, unions of concurrent lines, i.e., pencils and their small perturbations yield the extremal examples in several ranges where $t$ is small. In the Desarguesian case, recent stability results imply that these constructions are essentially unique whenever the number of atypical secants is sufficiently small.
	
	Second, when $q$ is a square, Baer subplanes provide examples with substantially fewer atypical secants than pencils for certain residue classes. These constructions explain the appearance of the $\sqrt q$-scale in our results and show that the general {lower} bound obtained from counting cannot be asymptotically improved in {several} parameter ranges. Further special point-line arrangements also emerge, in connection with so-called $(K;N)$-arcs; in our case $n$ takes the role of $K$ being the size of the point set, while $N$ denotes   the prescribed/maximum line-intersection number.
	
	Finally, when $n=qr-t$, for small positive $t$, optimal constructions appear to arise from unions of lines after deleting suitably chosen intersection points. This naturally leads to a new geometric extremal problem concerning line arrangements whose intersection points are themselves in general position. We formulate this notion as \emph{general position on level~2}. Somewhat unexpectedly, a natural restricted version of the problem leads to a Sidon-type problem over finite fields, thereby providing an arithmetic model for highly concentrated secant-size distributions.

	The problem above and our approach to it are closely related to various themes in combinatorial geometry and finite geometry. 
	The characterization of sets with no skew lines or tangents, and, more generally, sets without small secants is the subject of the study of covering and ($k$-fold) blocking sets, see for instance Bruen \cite{Br}, Blokhuis \cite{Blokhuis94},    the survey of Ball \cite{Ball03}, and the book of Kiss and Szőnyi \cite{Kiss-Szonyi}. Point sets having only a few line-intersection numbers form an active topic in finite geometry \cite{Durante} which is closely related to block designs, and notable examples are  subplanes, complete arcs, $(k;n)$-arcs and Korchmáros-Mazzocca arcs \cite{Kiss-Szonyi}.
	These are rigid configurations, the line-intersection pattern is supported on
	very few values. These examples suggest that a very small value of \(M(n,\Pi_q)\) should usually come
	from a highly organized geometric construction. One of the themes of the present paper is that this
	expectation is correct in several ranges of the size of $S$.  In particular, such conditions often invoke stability-type arguments, see the work of Szőnyi and Weiner \cite{Szonyi-Weiner}. 
	Avoiding a prescribed secant size can also be challenging. We refer to the theory of untouchable (or tangent-free) sets \cite{ Blokhuis, BSzW} or the recent paper of Héger and Nagy \cite{Heger-Nagy} for the general case.\\

	Our first general result is a universal {lower} bound obtained by the standard equation method, relying on double countings in the projective plane, and a refined application of the Cauchy--Schwarz inequality. Although this estimate is generally not asymptotically sharp, it has the correct order of magnitude for a broad  parameter range and closely reflects the behaviour exhibited by the geometric constructions. 
	
	\begin{theorem}\label{Cauchy-S1} Let  $n=qr+t$,  $0<n\le(q^2+q+1)/2$, where  $-q/2<t\le q/2.$  Then
		
		$$ \frac{(t(q+1)-r)^2}{qn-qr^2+(r-t)^2}\le M(n,\Pi_q).$$
	\end{theorem}

	{Similar statement hold for $n>(q^2+q+1)/2$ for complement symmetry reason (cf. Lemma 2.7). As a consequence, the Cauchy--Schwarz bound gives many values of $n$
		for which the number of atypical lines is large.}
	
	
	\begin{theorem}\label{CS-corr}
		For every \(\alpha\in[1,2]\), there are \(\Omega(q^{\alpha})\) values of \(n\) such that
		\[
		M(n,\Pi_q)=\Omega(q^{3-\alpha})
		.\]
	\end{theorem}
	
	

	Theorem \ref{Cauchy-S1} need not be tight in general, as we will show another natural lower bound on $M(n,\Pi_q)$, uniform in almost all $n$.
	
	\begin{proposition}
		For $n\in \{0,q+1, q^2,  q^2+q+1\}$, we have $M(n,\Pi_q)\in \{0,1\}$.
		Assume $1 \le n \le q^2+q$ and $n\neq q+1, q^2$. Then:
		\begin{itemize}
			\item[(1)] $M(n,\Pi_q) = q+1$, if $n \equiv 0, 1 \pmod q$;
			\item[(2)] $M(n,\Pi_q) > q+1$, if $n \not \equiv 0, 1 \pmod q$.
		\end{itemize} 
	\end{proposition}
	
	\begin{figure}
		\centering
		\includegraphics[width=\linewidth]{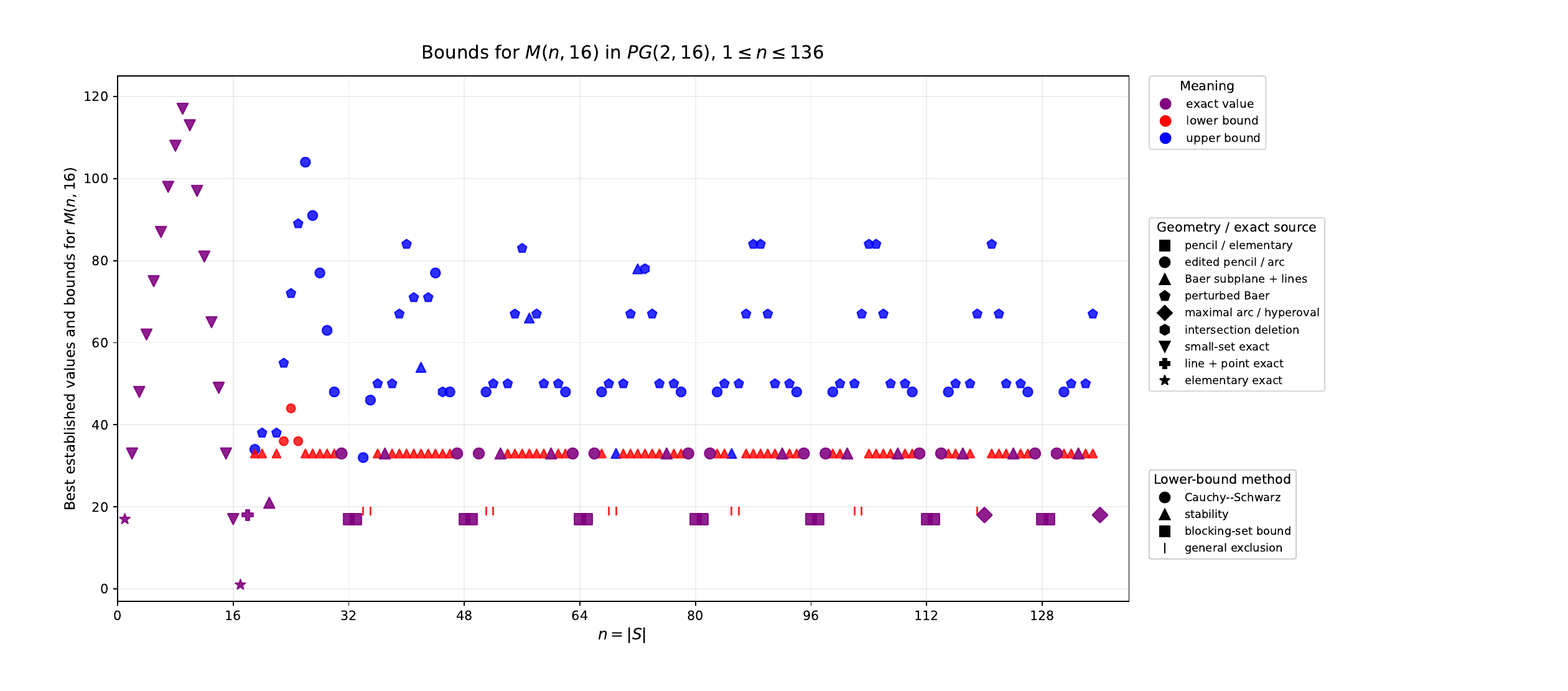}
		\caption{Lower and upper bounds for $q=16$ in the Desarguesian case, for $n<(q^2+q+1)/2$}
		\label{fig:placeholder}
	\end{figure}
	
	In the Desarguesian case, we will use much stronger stability-type results, which will enable us to characterize as well the sets where $M(n,\Pi_q)$ is small. 
	

	The paper is organized as follows.
	In Section~\ref{sec:preliminaries}, we collect notation, basic identities, and the external geometric
	results used later as black boxes. In Section~\ref{sec:general-bounds}, we discuss how small $M(n,\Pi_q)$ can  be and prove general lower bounds by counting and structural arguments, including the Cauchy--Schwarz bound and the
	blocking-set method. Then in 
	Section~\ref{sec:const}, we examine 
	which constructions yield matching or competitive upper bounds, and present explicit constructions
	and the resulting upper bounds on \(M(n,\Pi_q)\).  We  continue by studying under what conditions  these constructions are the only extremal configurations in Section~\ref{sec:desarguesian}. We examine  the Desarguesian case in detail, where earlier
	stability theorems allow stronger conclusions. 
	Finally, Section~\ref{sec:sidon} revisits the constructions for point sets with $n=qr+t$ with $-q/2<t<0$. Here a natural candidate construction, taking unions of lines and deleting suitably arranged intersection points, leads naturally to the concept of the so-called 
	general position on level 2 and to an additive combinatorial problem  related to  Sidon sets, which might be of independent interest.
	
	\section{Definitions and basic observations}\label{sec:preliminaries}
	
	Recall (Definition  \ref{def:k-secants} and \ref{def: s_m}) that $L_k=\{\ell\in \cL : n_\ell=k\}$ denotes the set of $k$-secants and $s_m(S)=\max_k \{|L_k(S)|\}$ denotes the maximum frequency of the secant sizes.
	\begin{definition}[Dominant lines and atypical lines] Let $S\subseteq \cP$. Choose a secant size $\kappa$ and call it a \textit{dominant secant size} 
		if
		$|L_\kappa(S)|=s_m(S)$, and call other secant sizes \emph{atypical}. A line $\ell$ is called dominant if $n_\ell=\kappa$
		for the chosen dominant secant size $\kappa$, and atypical otherwise.
	\end{definition}

	Note that the dominant secant size need not be unique in general, it is unique for every extremal point set
	when \(q\ge4\); see Proposition~\ref{prop:unique-dominant}.
	
	We recall our main extremal quantity from the introduction: $M(n,\Pi_q)$ is the minimum possible number of atypical lines among \(n\)-point sets.
	The investigation of the size of the set of atypical lines will capture better the dependence on $n$ and $q$.  
	
	
	\begin{definition} \
		\begin{itemize}
			\item  A pencil is a set of concurrent lines (often identified with their point set). A punctured  pencil is obtained from the point set of a pencil by removing the intersection point of its lines. A \em{carrier} is the intersection point of the lines of the pencil.
			\item A set of lines $\mathcal B$ in $\Pi_q$ is a \emph{covering line set} if every point of $\Pi_q$
			lies on at least one line of $\mathcal B$.
			{ Dually, a point set $B$ is a \emph{blocking set} if every line meets $B$.}
			\item A point set $\mathcal U$ is untouchable, or is a set without tangents if $|\ell\cap \mathcal U|\neq 1$ for every line $\ell$. 
			\item An \emph{oval} in a projective plane of order \(q\) is a set of \(q+1\) points no three of which are collinear.
			If \(q\) is even, a \emph{hyperoval} is a set of \(q+2\) points no three of which are collinear.
		\end{itemize}
	\end{definition}
	
	\begin{theorem}[Lower bounding the size of minimal blocking sets, Bruen-Pelikán, Blokhuis, \cite{Br, Blokhuis94}]\label{thm:blocking}
		Let $B$ be a blocking set in a projective plane of order $q$.
		If $B$ does not contain a line, then $
		|B|\ge q+\sqrt q+1.$
		Moreover, equality implies that $q$ is a square and $B$ is a Baer subplane.
		
		If $q$ is prime and the plane is Desarguesian, $\Pi_q=\PG(2,q)$, then every blocking set not containing a line has size at least $\frac{3(q+1)}2.$
	\end{theorem}
	
	\begin{theorem}[Blokhuis--Seress--Wilbrink]\label{thm:bsw}
		If $q$ is odd, then every nonempty set without tangents in $\PG(2,q)$ has size at least $
		q+\frac{\sqrt{2q}}{4}+2.$
	\end{theorem}
	
	As we indicated in the first section, the problem turns out to be connected to additive combinatorial concepts, most notably to Sidon sets.
	
	\begin{definition}
		Let $G$ be an Abelian  group and let $A\subseteq G$.  
		
		A set $A$ is called a \emph{Sidon set} (or \emph{$B_2$ set}) if the equation 
		$a+b=x$  ($a, b\in A$) has at most one unordered solution $(a,b)$ for all $x\in G$.
		
		More generally, a set $A$ is called a \emph{$B_2(h)$ set}
		(or \emph{$h$-Sidon set}) 
		if every \(x\in G\) has at most \(h\) unordered representations of the form \(x=a+b\) with \(a,b\in A\). 
		
		Finally,   we call \(A\) a \(B_2^*(h)\)-set (or \emph{$h$*-Sidon set}) if every \(x\in G\) has at most \(h\) unordered representations
		$$ x=a+b,\qquad a,b\in A,\quad a\ne b. $$
	\end{definition}
	
	\begin{result}\cite{Cille10}
		Let $G$ be a finite Abelian group of order $N$.
		If $A\subseteq G$ is a $B_2(h)$ set, then
		\[
		|A| \le \sqrt{2hN}+O(1),
		\] moreover, $|A| \le \sqrt{N}+1$ for $h=1$.
		Moreover, for every fixed $h$,  
		$$\max \{|A|: A\subseteq \mathbb Z_N \text { is } B_2(h)\}=\Omega(\sqrt{hN}).$$
		The same holds for $B_2^*(h)$-sets, as $$\max \{|A|: A\subseteq \mathbb Z_N \text { is } B_2^*(h)\}=\Omega(\sqrt{hN}).$$ \qed
	\end{result}
	
	\noindent To distinguish between the general case and the classical setting, we denote by $M^D(n, q)$ the minimum number of non-$\kappa$-secants specifically when the projective plane is Desarguesian.    

	We start with the simple complement symmetry.

	\begin{lemma}
		\label{lem: trivial M}\label{lem:complementer}
		$M(n,\Pi_q) = M(q^2+q+1-n, \Pi_q)$. \qed
	\end{lemma}
	\begin{proof}
		Observe that $L_{k}(S) = L_{q+1-k}(S^c)$, where $S^c$ is the complement of the set $S$ in the projective plane.  
		The claim follows.   
	\end{proof}
	
	Next we record the standard equalities.
	
	\begin{lemma}[Standard equalities]\label{lem:standard}
		Let $S\subseteq \Pi_q$ with $|S|=n$. Then
		\begin{align}
			\sum_k |L_k(S)| &= N,\label{eq:sum0}\\
			\sum_k k\,|L_k(S)| &= n(q+1),\label{eq:sum1}\\
			\sum_k k(k-1)\,|L_k(S)| &= n(n-1),\label{eq:sum2}\\
			\sum_k k^2\,|L_k(S)| &= n(n+q).\label{eq:sum3}
		\end{align}
	\end{lemma}

	\begin{lemma}\label{lem: n=0, 1}
		The following values hold for $M(n,\Pi_q)$ in any projective plane.
		\begin{itemize}
			\item[(1)] If $n = 0, q^2+q+1$ then $M(n,\Pi_q) = 0$; 
			\item[(2)] If $n=1$ or $n=q^2+q$, then $M(n,\Pi_q) = q+1$.
			\item[(3)] If $n \in \{q + 1, q^2\}$ then $M(n,\Pi_q)= 1$. 
		\end{itemize}
	\end{lemma}
	
	\begin{proof}
		$(1)$ and $(2)$ follow immediately from the definition of $M(n,\Pi_q)$. Counting the number of points in a set on the lines through a point in the set and a point not in the set shows that $M(n,\Pi_q)=0$ iff $n = 0, q^2+q+1$. The number of non-1-secants of the set of points of a line is $1$ (the line itself), which then proves $(3)$.     
	\end{proof}

	\section{General results and lower bounds for Problem \ref{maxi}}\label{sec:general-bounds}

	\subsection{A general Cauchy-Schwarz bound}
	
	Now we turn to the main general counting inequality. Recall that $N=(q^2+q+1)$.

	\begin{proposition}\label{Cauchy-S} Suppose that $S\subseteq \cP$ with $|S|=n=qr+t$, $0<n\le(q^2+q+1)/2$, where  $-q/2<t\le q/2.$  Then
		$$ s_m(S)\le \frac{n(n+q)N-n^2(q+1)^2}{Nr^2+n(q+n)-2n(q+1)r}$$
		
		i.e., 
		
		$$ M(n,\Pi_q)\ge \frac{(t(q+1)-r)^2}{qn-qr^2+(r-t)^2}.$$
	\end{proposition}

	\begin{proof}
		Let us choose \(\kappa\) such that
		$ |\mathcal L_\kappa(S)|=s_m(S).$ Put
		
		$$ m:=N-|\mathcal L_\kappa(S)|=N-s_m(S). $$
		
		Thus exactly \(m\) lines lie outside \(\mathcal L_\kappa(S)\).
		Writing up the standard equations, cf. Lemma \ref{lem:standard},
		we obtain that

		\begin{equation}
			\sum_{\ell\not\in L_{\kappa}} n_{\ell} = n(q+1)-(N-m)\kappa, \ \
			\sum_{\ell\not\in L_{\kappa}} n_{\ell}(n_{\ell}-1) = n(n-1)-(N-m)\kappa(\kappa-1).
		\end{equation}
		
		The Cauchy-Schwarz inequality then implies

		\begin{equation}
			\begin{split}
				(\sum_{\ell\not\in L_{\kappa}} n_{\ell})^2 &\leq \sum_{\ell\not\in L_{\kappa}} n_{\ell}^2\cdot m\\&= [(n(q+1)-(N-m)\kappa)+(n(n-1)-(N-m)\kappa(\kappa-1))]\cdot m\\
				& =[n(q+n)-(N-m)\kappa^2)]\cdot m.
			\end{split}
		\end{equation}
		Note that in the last step, we used the fact that $\sum_{\ell\not\in L_{\kappa}} n_{\ell}^2 = \sum_{\ell\not\in L_{\kappa}} n_{\ell}(n_{\ell}-1)+\sum_{\ell\not\in L_{\kappa}} n_{\ell}$ clearly holds.
		Hence,
		$$(n(q+1)-(N-m)\kappa)^2\le [n(q+n)-(N-m)\kappa^2)]\cdot m, $$ that is,

		$$  (N-m)(N\kappa^2+n(q+n)-2n(q+1)\kappa) \le (-n^2(q+1)^2+n(n+q)N).$$
		
		\noindent The minimum of $m$ could be derived from here by minimizing the coefficient  of $(N-m)$ $$(N\kappa^2+n(q+n)-2n(q+1)\kappa)$$ in $\kappa$. Since this is a quadratic inequality, the optimal choice can be derived studying the corresponding discriminant. Since $n<N/2$ and $n=qr+t$ with $-q/2<t\le q/2$,
		we have
		\[
		\left|\frac{n(q+1)}{N}-r\right|<\frac12.
		\]
		Thus $r$ is the integer minimizing
		$A(k):=Nk^2+n(n+q)-2n(q+1)k.$
		In particular, $A(\kappa)\ge A(r)$ for every
		dominant secant size $\kappa$. Substituting this
		into the preceding inequality yields the claim.
	\end{proof}
	\begin{remark}
		Equality in the upper bound  of Proposition \ref{Cauchy-S} on $s_m$ 
		can occur only if all atypical lines have the same intersection size with \(S\). Thus equality in the resulting lower bound for \(M(n,\Pi_q)\) requires, in particular, that the dominant secant size be the minimizing integer \(r\) and that $ n_\ell $ be constant over all \(\ell\notin\mathcal L_r(S)\). Such an equality case therefore gives rise to a highly regular incidence structure, closely related to a block design.
	\end{remark}
	\begin{proof}[Proof of Theorem \ref{CS-corr}]
		{Let
			\[
			R:=\max\{1,\lfloor c q^{\alpha-1}\rfloor\},
			\]
			where \(c>0\) is a sufficiently small constant, and choose
			\(1\le r\le R\) and \(\lceil q/8\rceil\le t\le \lfloor q/4\rfloor\)}
		{ and hence} 
		$
		t=\Omega (q).
		$
		Proposition \ref{Cauchy-S} ensures that for
		$
		n=qr+t
		$
		we have
		$$
		M(n,\Pi_q)\ge
		\frac{(t(q+1)-r)^2}{qn-qr^2+(r-t)^2}.
		$$
		
		We estimate the numerator and denominator separately. Since \(t\ge q/8\) and \(r\le R=O( q^{\alpha-1})\)
		, we have
		$
		t(q+1)-r\gg q^2,
		$
		and therefore
		$
		(t(q+1)-r)^2\gg q^4.
		$
		
		On the other hand,
		$
		qn-qr^2+(r-t)^2
		=
		q\,r(q-r)+qt+(r-t)^2.
		$
		As \(r\le R\), we have
		$
		q\,r(q-r)\le q^2r\ll q^{\alpha+1},
		$
		while
		$
		qt\ll q^2,$ and $
		(r-t)^2\ll q^2
		$ holds. 
		Since \(\alpha\ge 1\), it follows that
		$
		qn-qr^2+(r-t)^2\ll q^{\alpha+1}.
		$
		Consequently,
		\[
		M(n,\Pi_q)\gg \frac{q^4}{q^{\alpha+1}} = q^{3-\alpha}.
		\]
		
		It remains to count the number of admissible values of \(n\). There are \(R=\Omega( q^{\alpha-1})\) admissible values for \(r\) and \(\Omega( q)\) possible choices for \(t\), hence \(\Omega(q^\alpha)\) pairs \((r,t)\). Since
		$
		n=qr+t
		$
		and \(0<t<q\), distinct pairs \((r,t)\) yield distinct values of \(n\). Thus there are \(\Omega(q^\alpha)\)
		values of \(n\) for which
		$
		M(n,\Pi_q)\gg q^{3-\alpha}.
		$
	\end{proof}
	
	\subsection{Stability results when the set has small size}
	
	\begin{lemma}
		\label{lem: S small} Suppose that $t$-arcs exist in $\Pi_q$ for all $t\le (2-\sqrt2)q+1$. (Note that for Desarguesian planes, this holds in a strong sense, for all  $t\le q+1)$.\\
		If $n \le \left\lfloor 2q+\frac{3}{2}- \sqrt{2q^2+4q+\tfrac14}\right\rfloor \approx(2-\sqrt 2)q$, then $M(n,\Pi_q) =  M(q^2+q+1-n, \Pi_q) = n(q+1)- {n \choose 2}$, while if  $q\ge n >\left\lfloor 2q+\frac{3}{2}- \sqrt{2q^2+4q+\tfrac14}\right\rfloor$, then $M(n,\Pi_q) = M(q^2+q+1-n, \Pi_q) =$ {$q^2+q+1$}$-qn$.  
	\end{lemma}
	
	\begin{proof}  The transition point will be determined   by comparing $|L_0|$ and $|L_1|$. Suppose first that $1<n\le \left\lfloor 2q+\frac{3}{2}- \sqrt{2q^2+4q+\tfrac14}\right\rfloor $ and $S$ is an arbitrary $n$-element set.
		Then    $|L_k(S)|\le \binom{n}{2}$ for all $k\ge 2$, while $|L_1(S)|\le nq$ holds. Let $S_0$ be an $n$-arc.  Hence it is enough to show that $|L_0(S_0)|=q^2+q+1-(q+1)n+\binom{n}{2}\ge \max \{nq, \binom{n}{2}\} $ holds. The inequality in turn follows from the upper bound  on $n.$ \\
		On the other hand,  suppose next that $q\ge n> \left\lfloor 2q+\frac{3}{2}- \sqrt{2q^2+4q+\tfrac14}\right\rfloor $ and $S$ is an arbitrary $n$-element set. Then    $|L_k(S)|\le \binom{n}{2}$ for all $k\ge 2$, while $|L_0(S)|\le q^2+q+1-(q+1)n+\binom{n}{2}$ holds. 
		Let $S_1$ be a set of collinear points. Hence it is enough to show that $|L_1(S_1)|=qn\ge \max \{q^2+q+1-(q+1)n+\binom{n}{2}, \binom{n}{2}\} $ holds. The inequality in turn follows from the upper bound  on $n.$
	\end{proof}
	
	\begin{lemma}
		\label{lem:not dual bl}
		If for a set $S$, the set of lines $\mathcal{L}({\Pi_q})\setminus L_{\kappa}(S)$ does not cover all points of the projective plane then either $|S| = (q+1)\kappa(S)$ or  $|S^c| = (q+1)(q+1 - \kappa(S))$. \qed 
	\end{lemma}

	\begin{proof} Let \(P\) be a point incident only with \(\kappa(S)\)-secants. If \(P\notin S\), then the \(q+1\) lines through \(P\) partition \(S\), whence
		$ |S|=(q+1)\kappa(S). $\\
		If \(P\in S\), then every line through \(P\) contains \(\kappa(S)-1\) further points of \(S\), and hence
		$ |S|=1+(q+1)(\kappa(S)-1). $
		Equivalently,
		$|S^c|=(q+1)(q+1 - \kappa(S)).$
	\end{proof}

	\begin{lemma}
		\label{lem:untouchable}
		Let $S$ be a non-empty set in $\Pi_q$ and assume that neither $S$ nor $S^c$ is a line. Suppose also that the non-$\kappa(S)$-secants do not cover all points of the plane. Then through each point there pass $0$ or at least $2$ non-$\kappa(S)$-secants of $S$. In particular, $S$ has either $0$ or at least $q+2$ non-$\kappa(S)$-secants.
	\end{lemma}
	
	\begin{proof}
		By Lemma \ref{lem:not dual bl} and Lemma \ref{lem:complementer} we may suppose that $|S|= \kappa(S)(q+1)$ (i.e., there exists a point outside $S$ with only $\kappa(S)$-secants through it). If now $Q$ is a point with exactly one non-$\kappa(S)$-secant line $\ell$ through it, then $Q$ must be in $S$. Counting the points of $S$ on the lines through $Q$, we get that $|\ell \cap S| = q+ \kappa(S)$ and so $\kappa(S) = 0, 1$. Since $S$ is non-empty with size $|S|= \kappa(S)(q+1)$, $\kappa(S)$ must be $1$ and $|S|=q+1$. Again looking around $Q$ we see that $S$ is a line, which contradicts our assumption.
		Hence through each point there pass $0$ or at least $2$ non-$\kappa(S)$-secants of $S$. If a non-$\kappa(S)$-secant exists then counting the non-$\kappa(S)$-secants through the points of this secant we see at least $q+2$ of them.
	\end{proof}
	
	\begin{lemma}
		\label{lem:pencil}
		Let $S$ be a set  with less than $2q + 1$ lines not in $L_{\kappa(S)}$. Assume also that among these lines there are $q+1$ through a point $P$. Then $S$ is:
		\begin{itemize}
			\item[(1)] a pencil or a punctured pencil,
			\item[(2)] a point set of a line and a disjoint point 
			\item[(3)] the complement of $(2)$.
		\end{itemize}
	\end{lemma}
	
	\begin{proof}
		Note that on each line through $P$ there must exist a point $Q \neq P$ such that $PQ$ is the unique non-$\kappa(S)$-secant through $Q$. Assume that $Q_1$ and $Q_2$ are such points so that the line $PQ_1$ is different from $PQ_2$. If $Q_1 \in S$ but $Q_2 \not \in S$, then counting the number of points of $S$ on the lines through $Q_1$ and through $Q_2$, we get that $$|S|=\kappa q+|PQ_2\cap S|=(\kappa-1)q+|PQ_1\cap S|,$$ thus $|S\cap \{PQ_1 \setminus P \}| = q$ and  $|S\cap \{PQ_2 \setminus P \}| = 0$; and hence $S$ is a pencil or a punctured pencil.
		
		Otherwise we see the same number of points of $S$ on each line through $P$.  Write
		$a=|(\ell\setminus\{P\})\cap S|$.
		The cases $a=0,q$ again give pencils or punctured pencils,
		so assume $0<a<q$. By complementation, we may assume
		$P\notin S$. Then $n=a(q+1)$.  As mentioned before on each line through $P$, there exists a point $R$, through which $RP$ is the only non-$\kappa$-secant. Clearly, $R \in S$ otherwise we get that $\kappa(S) = a$ and hence the lines through $P$ are in $L_{\kappa(S)}$. Counting the points of $S$ on the lines through $R$ we get that $\kappa(S)-1 = a$. 
		
		Set $S'=S\cup\{P\}$. All lines through $P$ are now
		$\kappa$-secants, and fewer than $q$ other lines are not.
		In particular, $\kappa$ is dominant for $S'$.
		These exceptional lines do not cover the plane.
		By Lemma~\ref{lem:untouchable} and
		Lemma~\ref{lem: n=0, 1}, $S'$ is a line, the complement
		of a line, the empty set, or the whole plane.
		Since $|S'|=a(q+1)+1$ and $0<a<q$, only the complement
		of a line is possible.
		Thus $S$ is the complement of a line together with
		an external point. Undoing complementation gives
		the remaining possibility. \end{proof}


	\begin{lemma} \label{lem:0,1uj}
		Assume $1 < n < q^2+q$ and $n\neq q+1, q^2$, then:
		\begin{itemize}
			\item[(1)] $M(n,\Pi_q) = q+1$, if $n \equiv 0, 1 \pmod q$;
			\item[(2)] $M(n,\Pi_q) > q+1$, if $n \not \equiv 0, 1 \pmod q$.
		\end{itemize} 
	\end{lemma}
	
	\begin{proof}
		Let $S$ be a set of size $n$. If its non-$\kappa(S)$-secants
		do not cover the plane, then Lemma \ref{lem:untouchable}, Lemma \ref{lem: n=0, 1} and our assumption on $n$ give at least $q+2$ such lines. If they cover the plane, then by Theorem~\ref{thm:blocking} there are at least $q+1$ of them, with equality only if they form a pencil. Thus $M(n,\Pi_q)\ge q+1$.
		
		If equality holds, the non-$\kappa(S)$-secants form a pencil, so
		Lemma~\ref{lem:pencil} applies. A line plus an external
		point and its complement have $q+2$ non-dominant secants;
		hence $S$ is a pencil or a punctured pencil, and
		$n\equiv0,1\pmod q$. If $n\equiv0,1\pmod q$, we can choose a pencil or a
		punctured pencil of size $n$ with exactly $q+1$
		non-$\kappa(S)$-secants.
	\end{proof}

	\begin{lemma}
		\label{lem: hyperoval}
		Assume that a set $S$ admits
		$q+2$ non-$\kappa(S)$-secants, $l_1, l_2, \dots l_{q+2}$ forming a dual hyperoval. Then $S$ is one of the following
		\begin{itemize}
			\item[(1)] the set of all points lying on one of the lines $l_i$, $1 \le i \le q+2$;
			\item[(2)] the set of all points lying on one of the lines $l_i$, $1 \le i \le q+2$, except those on a fixed line $l_j$.
			\item [(3)] the complement of either $(1)$ or $(2)$.
		\end{itemize}    
	\end{lemma}
	
	\begin{proof}
		For any point not lying on any of the lines $l_i$, we see only $\kappa(S)$-secants of $S$ through it and so either all points outside the lines $l_i$ are in $S$ or none. By taking the complement if necessary, we may assume that none of these points belongs to $S$. Counting through any such point gives $|S|=(q+1)\kappa(S)$.
		
		If no $l_i$ contains both a point of $S$ and a point
		outside $S$, then each $l_i$ is either contained in
		$S$ or disjoint from it. Since the $l_i$ meet pairwise,
		they must all be of the same kind. They cannot all be
		disjoint, as then $S$ would be empty. Thus $(1)$ holds. Now pick a line $l_i$ which is not a $\kappa(S)$-secant of $S$ and assume that there exist points $Q \in l_i \cap S$ and $R \in \{l_i \setminus S\}$. Let $l_Q$ and $l_R$ be the other non-$\kappa(S)$-secants of $S$ through $Q$ and $R$. Counting the points of $S$ on the lines through $Q$ and through $R$, we get that either $|l_Q \cap S| = q+ 1$ and $|l_R \cap S| = 1$ or  $|l_Q \cap S| = q$ and $|l_R \cap S| = 0$. 
		
		More generally, counting through the intersection of
		two distinct lines $l_a,l_b$ among the $l_i$ gives
		\[
		|l_a\cap S|+|l_b\cap S|=
		\begin{cases}
			2\kappa(S)+q,& l_a\cap l_b\in S,\\
			2\kappa(S),& l_a\cap l_b\notin S.
		\end{cases}
		\]
		
		Suppose first that $|l_Q\cap S|=q$ and
		$|l_R\cap S|=0$. Since $l_R$ is disjoint from $S$,
		the displayed equality gives
		$|l_j\cap S|=2\kappa(S)$ for every $l_j\ne l_R$.
		In particular, $2\kappa(S)=|l_Q\cap S|=q$.
		Thus every $l_j\ne l_R$ contains all its points outside
		$l_R$ in $S$, giving (2).
		
		In the other case, $|l_Q\cap S|=q+1$ and
		$|l_R\cap S|=1$. The point $l_Q\cap l_R$ belongs to $S$,
		so the displayed equality gives $\kappa(S)=1$.
		For any $l_j\ne l_Q,l_R$, the point $l_j\cap l_R$
		is not in $S$, and hence $|l_j\cap S|=1=\kappa(S)$.
		This contradicts the choice of the $l_j$.
		Taking complements gives $(3)$.
	\end{proof}

	For $q=2$, Lemmas~\ref{lem: n=0, 1} and
	\ref{lem:0,1uj} already determine all values of $M(n,\Pi_2)$. For $q\ge3$, The next corollary follows from the preceding lemmas and the dual form of Theorem~\ref{thm:blocking}

	\begin{corollary}\label{cor: gen M} Suppose that $q\ge3$ and $0<n<q^2+q+1$. Then
		$M(n,\Pi_q)>q+2$, unless $n\equiv0,1\pmod q$,
		$n\in\{q+2,q^2-1\}$, or $\Pi_q$ admits a dual
		hyperoval and
		\[
		n\in\left\{
		\frac{(q+1)(q+2)}2,\,
		\frac{q(q+1)}2,\,
		\frac{q(q-1)}2,\,
		\frac{q^2+q+2}2
		\right\}.
		\]
	\end{corollary}

	\subsection{Structural arguments  relying on blocking sets}

	\begin{proposition} \label{covering}
		
		Let \(S\subseteq\Pi_q\), and let \(\kappa\) be a dominant secant size. If
		
		$$ |S|\not\equiv0,1\pmod{q+1}, $$
		
		then the non-\(\kappa\)-secants form a covering line set.
	\end{proposition}
	\begin{proof}
		
		If some point is not covered by a non-\(\kappa\)-secant, Lemma \ref{lem:not dual bl} gives either
		$|S|=(q+1)\kappa $
		or
		$ |S^c|=(q+1)(q+1 -\kappa(S)). $
		In the first case \(|S|\equiv0\pmod{q+1}\), while in the second
		$|S|=q^2+q+1-(q+1)(q+1 -\kappa(S))\equiv1\pmod{q+1}. $
		This proves the claim.
	\end{proof}

	\begin{proposition} \label{covering1} Let \(S\subseteq\Pi_q\), let \(\kappa\) be a dominant secant size of \(S\), and suppose that
		
		$$ |S|\not\equiv0,1\pmod{q+1}. $$
		
		If \(S\) is neither a pencil nor a punctured pencil, then
		
		$$ N-|\mathcal L_\kappa(S)|\ge q+\sqrt q+1. $$
		
	\end{proposition} 
	\begin{proof}
		Put
		$ \mathcal B=\mathcal L\setminus\mathcal L_\kappa(S). $
		By Proposition \ref{covering}, \(\mathcal B\) is a covering line set. Suppose
		that 
		$ |\mathcal B|<q+\sqrt q+1. $
		Choose an inclusion-minimal covering subset \(\mathcal B_0\subseteq\mathcal B\). By dual Bruen, either \(|\mathcal B_0|\ge q+\sqrt q+1\), contradiction, or \(\mathcal B_0\) contains a full pencil. But because \(\mathcal B_0\) is minimal, in the latter case it is exactly that pencil. Hence \(\mathcal B\) contains all \(q+1\) lines through some point. Since
		$ |\mathcal B|<q+\sqrt q+1<2q+1, $
		Lemma \ref{lem:pencil} applies and says that \(S\) is a pencil, a punctured pencil, a line plus an external point, or the complement of the latter. The first two are excluded; the last two have cardinalities \(q+2\equiv1\pmod{q+1}\) and \(q^2-1\equiv0\pmod{q+1}\), so they are also excluded.
	\end{proof}
	\begin{proposition} Suppose that 
		$n \not \equiv 0,1 \pmod q, n\not \equiv0,1 \pmod {q+1}$ then $M(n,\Pi_q)\ge q+\sqrt{q}+1.$
	\end{proposition}
	\begin{proof}
		Choose an \(n\)-point set \(S\) attaining \(M(n,\Pi_q)\). A pencil has size \(rq+1\), while a punctured pencil has size \(rq\), so \(n\not\equiv0,1\pmod q\) excludes both possibilities. Since also \(n\not\equiv0,1\pmod{q+1}\), Propositions \ref{covering} and \ref{covering1} apply and yield the result.
	\end{proof}

	\section{Upper bounds on $M(n,\Pi_q)$: constructions}\label{sec:const}
	
	In view of Lemma~\ref{lem: n=0, 1}, we may assume throughout this section that
	$n\notin\{0,q+1,q^2,q^2+q+1\}$.
	
	\begin{construction}\label{constrs}
		We consider the following point sets:
		\begin{enumerate}[(i)]
			\item a punctured pencil;
			\item a pencil;
			\item a punctured pencil together with $t>0$ points of an arc containing its carrier;
			\item the union of $r_0\ge2$ lines through a common point and $r_1=r-r_0\ge1$ further lines, with all pairwise intersections distinct apart from the common point of the first $r_0$ lines, after deleting a subset of the intersection points;
			\item a maximal $((q+1)k-q;k)$-arc (a Denniston arc), when $q=2^h$, the plane is Desarguesian, and $k\mid q$~\cite{Kiss-Szonyi};
			\item the union of $k$ pairwise disjoint Baer subplanes and $s$ concurrent lines, as specified in Proposition~\ref{Baer+line};
			\item a perturbation of a set in~(vi), obtained by adding or deleting fewer than $\sqrt q$ points;
			\item the union of $r>1$ lines in general position, with a subset of their intersection points deleted. In the level-two case, no line outside the original $r$ lines contains three of their pairwise intersection points.
		\end{enumerate}
	\end{construction}
	
	The first three constructions, together with maximal arcs, yield the standard bounds below. The partial-pencil construction in~(iv) is treated explicitly after Proposition~\ref{upperb}; the Baer-subplane construction in~(vi) is analyzed in Proposition~\ref{Baer+line}. Construction~\ref{constrs}(viii) will lead to Lemma~\ref{gen_pos} and the level-two problem of Section~\ref{sec:sidon}.
	
	\begin{proposition}[Upper bounds on $M(n,\Pi_q)$]\label{upperb} \ 
		\begin{enumerate}[(i)]
			\item For $1\le r\le q+1$,
			$ M(qr,\Pi_q)\le q+1.$
			\item For $0\le r\le q$,
			$ M(qr+1,\Pi_q)\le q+1.$
			\item If $2\le t\le q/2$ and $1\le r\le q+2-t$, and a $t$-arc exists in $\Pi_q$, then
			\[
			M(qr+t,\Pi_q)\le tq+1-\binom{t-1}{2}.
			\]
			\item Suppose that $1\le m\le r-2$ and that $r$ lines $\ell_1,\ldots,\ell_r$ have the following incidence pattern: $\ell_1,\ldots,\ell_{r-m}$ pass through a common point, and all their other pairwise intersections are distinct. Write
			$ X=\{\ell_i\cap\ell_j:1\le i<j\le r\},$ $ 
			x_m=|X|=1+{m(2r-m-1)}/2.$
			For every integer $t$ satisfying
			$ x_m-m-2\le t\le 2x_m-m-2$,
			we have the  bound
			\[
			M(rq-t,\Pi_q)\le x_m(q-1)+m+2.
			\]
			If, in addition, no line other than $\ell_1,\ldots,\ell_r$ contains more than two points of $X$, this improves to
			\[
			M(rq-t,\Pi_q)\le
			x_m(q+1)-\binom{x_m}{2}
			+(r-m)\binom m2+m\binom{r-2}{2}.
			\]
			\item If $q=2^h$, $\Pi_q$ is Desarguesian, and $r+1\mid q$, then
			\[
			M(qr+r+1,\Pi_q)\le q\left(\frac{q+1}{r+1}-1\right).
			\]
			\item Under the hypotheses of Proposition~\ref{Baer+line}, put
			$ n=(k+s)q+k(\sqrt q-s)+(k+s).$
			Then
			\[
			M(n,\Pi_q)\le
			\bigl((k-1)(s+1)+2\bigr)q
			-(k-1)(s-1)\sqrt q+k.
			\]
		\end{enumerate}
	\end{proposition}
	
	\begin{proof}[Proof of Proposition~\ref{upperb}(iii)]
		Choose a $t$-arc containing the carrier $P$ of a punctured pencil of $r$ lines, none of which contains another point of the arc; this is possible since $r\le q+2-t$. After adding the $t$ arc points, only the $q+1$ lines through $P$ and the lines through at least one of the other $t-1$ arc points can cease to be $r$-secants. Among the lines not through $P$, their number is
		\[
		(t-1)q-\binom{t-1}{2},
		\]
		since no three of the added points are collinear. The claim follows.
	\end{proof}
	
	\begin{proof}[Proof of Proposition~\ref{upperb}(iv)]
		The union of the $r$ lines has
		\[
		r(q+1)-(r-m-1)-(x_m-1)=rq+m+2-x_m
		\]
		points. Deleting
		\[
		d=t-(x_m-m-2),\qquad 0\le d\le x_m,
		\]
		points of $X$ gives a set $S$ with $|S|=rq-t$.
		
		Every line avoiding $X$ meets the $r$ original lines in $r$ distinct points, none deleted from $S$. Thus
		\[
		M(rq-t,\Pi_q)\le
		|\{\ell:\ell\cap X\ne\varnothing\}|.
		\]
		To obtain a bound without further hypotheses, count the $r$ original lines, the at most $q+1-(r-m)$ further lines through their common point, and at most $q-1$ additional lines through each of the remaining $x_m-1$ points of $X$. Hence
		\[
		|\{\ell:\ell\cap X\ne\varnothing\}|
		\le r+q+1-(r-m)+(x_m-1)(q-1)
		=x_m(q-1)+m+2.
		\]
		
		Under the additional hypothesis, the only lines containing at least three points of $X$ are the original $\ell_i$. Counting the lines through $X$ by first summing the $x_m(q+1)$ point-line incidences and then subtracting one for each pair of points of $X$, we must restore $\binom{j-1}{2}$ for each line containing $j\ge3$ points of $X$. The first $r-m$ original lines contain $m+1$ points of $X$ each, and the remaining $m$ lines contain $r-1$ such points each. Consequently,
		\[
		|\{\ell:\ell\cap X\ne\varnothing\}|
		=x_m(q+1)-\binom{x_m}{2}
		+(r-m)\binom m2+m\binom{r-2}{2},
		\]
		as required.
	\end{proof}
	
	We next turn to the construction involving Baer subplanes and concurrent lines.
	
	\begin{proposition}\label{Baer+line}
		Suppose that $q$ is a square and $\B_1,\ldots,\B_k$ are $k>0$ pairwise disjoint Baer subplanes of $\PG(2,q)$. Let $\ell_1,\ldots,\ell_s$ be $s\ge0$ lines such that $|\ell_i\cap\B_j|=1$ for every $i,j$. If $s\ge1$, assume that these lines are concurrent at a point $P\in\B_1$. Put
		\[
		S=\bigcup_{j=1}^{k}\B_j\;\cup\bigcup_{i=1}^{s}\ell_i.
		\]
		Then
		\[
		|S|=(k+s)q+k(\sqrt q-s)+(k+s),
		\]
		and
		\[
		|L_{k+s}(S)|\ge N-
		\left[\bigl((k-1)(s+1)+2\bigr)q
		-(k-1)(s-1)\sqrt q+k\right].
		\]
	\end{proposition}
	
	\begin{proof}
		Since the Baer subplanes are disjoint,
		\[
		\left|\bigcup_{j=1}^{k}\B_j\right|=k(q+\sqrt q+1).
		\]
		Each of the $s$ lines adds exactly $q+1-k$ points outside the union of the subplanes; these added point sets are disjoint because the lines meet only at $P\in\B_1$. This proves the cardinality formula (and is immediate when $s=0$).
		
		Call a line \emph{fine} if it meets each $\B_j$ and each $\ell_i$ in exactly one point, with all these points distinct. Every fine line is a $(k+s)$-secant of $S$. If $s=0$, every non-fine line is a Baer-subline of one of the $\B_j$, so there are at most $k(q+\sqrt q+1)$ such lines, which is stronger than the stated bound.
		
		Suppose $s\ge1$. A non-fine line must be a Baer-subline of some $\B_j$, pass through $P$, or pass through $\ell_i\cap\B_j$ for some $i$ and $j\ge2$. There are at most $k(q+\sqrt q+1)$ lines of the first type. Apart from the Baer-sublines of $\B_1$, at most $q-\sqrt q$ lines pass through $P$. Similarly, for each of the $(k-1)s$ points $\ell_i\cap\B_j$ with $j\ge2$, at most $q-\sqrt q$ lines are not Baer-sublines of $\B_j$. By the union bound, the number of non-fine lines is therefore at most
		\[
		k(q+\sqrt q+1)+(q-\sqrt q)+(k-1)s(q-\sqrt q)
		=\bigl((k-1)(s+1)+2\bigr)q
		-(k-1)(s-1)\sqrt q+k.
		\]
		This yields the asserted bound on $|L_{k+s}(S)|$.
	\end{proof}
	
	In particular, if a set $S_0$ constructed above has size $n_0$ and at least $N-B$ $(k+s)$-secants, adding or deleting $d$ points changes only lines through those points, giving
	\[
	M(n_0\pm d,\Pi_q)\le B+d(q+1)
	\]
	whenever the indicated addition or deletion is possible. This is the perturbation principle in Construction~\ref{constrs}(vii).
	
	The next lemma concerns Construction~\ref{constrs}(viii). Its hypothesis is precisely general position on level two, as defined in Section~\ref{sec:sidon}.
	
	\begin{lemma}\label{gen_pos}
		Let $\ell_1,\ldots,\ell_r$ be $r\ge3$ lines in general position in $\PG(2,q)$, and put
		\[
		X=\{\ell_i\cap\ell_j:1\le i<j\le r\}.
		\]
		Assume that no line outside $\{\ell_1,\ldots,\ell_r\}$ contains three points of $X$. Then the number $b(X)$ of lines meeting $X$ is
		\[
		b(X)=r+\binom r2(q-1)-3\binom r4
		=r+\binom r2(q-1)
		-\frac12\binom{r-2}{2}\binom r2.
		\]
		Consequently, whenever
		\[
		\binom r2-r\le t\le r(r-2),
		\]
		we have
		\[
		M(rq-t,\Pi_q)\le b(X).
		\]
	\end{lemma}
	
	\begin{proof}
		The original $r$ lines meet $X$. Each of its $\binom r2$ points lies on a further $q-1$ lines. Under the hypothesis, an additional line meets at most two points of $X$; hence the only repeated counts correspond to pairs of intersection points associated with four distinct original lines. There are $3\binom r4$ such pairs, giving the formula for $b(X)$.
		
		The union of the $r$ lines has $r(q+1)-\binom r2$ points. Deleting any $d$ points of $X$, with $0\le d\le\binom r2$, gives a set of size $rq-t$, where
		\[
		t=\binom r2-r+d.
		\]
		Every line avoiding $X$ remains an $r$-secant, and the claimed bound follows.
	\end{proof}
	
	Finally, a simple pencil perturbation shows that extremal sets have a unique dominant secant size.
	
	\begin{proposition}\label{prop:unique-dominant}
		Let $\Pi_q$ be a projective plane of order $q\ge4$. Then every extremal point set has a unique dominant secant size.
	\end{proposition}
	
	\begin{proof}
		Put $N=q^2+q+1$. The empty set and the full point set clearly have a unique dominant secant size. By complementation, it therefore suffices to consider $0<n\le N/2$. Write $n=rq+t$, where $0\le t<q$.
		
		If $t=0$, a punctured pencil of $r$ lines has $q^2$ $r$-secants, and thus
		\[
		M(n,\Pi_q)\le q+1<\frac N2.
		\]
		Suppose $t\ge1$. Fix a point $P$ and $r$ lines through $P$, and put
		\[
		S_0=\{P\}\cup\bigcup_{i=1}^{r}(\ell_i\setminus\{P\}),
		\qquad |S_0|=rq+1.
		\]
		Every line not through $P$ is an $r$-secant of $S_0$. Adding $t-1$ further points destroys at most $q(t-1)$ of these $q^2$ secants. Therefore
		\[
		M(n,\Pi_q)\le qt+1.
		\]
		Applying the same bound to the complementary cardinality yields
		\[
		M(n,\Pi_q)\le q\min\{t,q+1-t\}+1. \tag{*}
		\]
		For even $q$ the right-hand side is at most $q^2/2+1<N/2$. For odd $q$, the same holds unless $t=(q+1)/2$.
		
		In the remaining case $q\ge5$ and $t-1\ge2$. Since $n\le N/2$, we have $r\le(q-1)/2$, so at least two lines through $P$ lie outside the chosen pencil. Choose two of the added points on distinct such lines. The line joining them does not pass through $P$ and is counted twice among the affected non-$P$ lines. Hence at most $q(t-1)-1$ of the former $r$-secants are destroyed, and
		\[
		M(n,\Pi_q)\le qt
		=\frac{q(q+1)}2=\frac{N-1}{2}<\frac N2.
		\]
		Thus $M(n,\Pi_q)<N/2$ for every $n$. If $S$ is extremal, then $s_m(S)=N-M(n,\Pi_q)>N/2$. Since distinct secant-size classes are disjoint, at most one of them can contain more than $N/2$ lines.
	\end{proof}
	
	\section{Erd\H os's two-coloring problem and bounds on $M(n,\Pi_q)$ and $M^D(n, q)$} \label{sec:desarguesian}
	
	In arbitrary projective planes the counting arguments of the previous sections already imply nontrivial lower bounds on $M(n,\Pi_q)$. In Desarguesian planes, however, considerably stronger conclusions become possible. Recent stability theorems classify point sets with few atypical secants, and therefore translate almost immediately into exact values and sharp lower bounds for $M(n,\Pi_q).$ These results originated from the coloring problem of Erdős.
	
	{
		\begin{problem}[Erdős] \label{prob: Erdos}
			What is the maximum number of lines in the projective plane $\Pi_q$
			that have the same number of red and blue points in two-colorings of the point set of $\Pi_q$? 
		\end{problem}
		This problem has been studied and generalized recently by Blokhuis, Markó, Szőnyi and Weiner \cite{BMSzW}, who {introduced the } concept of $(r, b)$-colored lines.  A line $\ell$ is $(r,b)$-colored in the  two-coloring of the plane  if the number of
		red points on $\ell$ is $r$, while the number of blue points is $b$. The natural question then is to determine the structure (and stability) of the two-colorings of   point sets {that } attain the maximum number of $(r,b)$-colored lines. If $S$ is the set of red points, the $(r,b)$-colored
		lines are precisely the $r$-secants of $S$. In the
		coloring problem $r$ is fixed, whereas the size of
		$S$ need not be fixed.
	}

	Stability results that characterize colorings with not too many non-$(r, b)$-colored lines provide direct exact values and lower bounds for $M(n,\Pi_q)$. Such characterizations identify the "candidate" configurations for a fixed $n$ and hence $M(n,\Pi_q)$ is found by minimizing the number of non-typical secants over these specific structures. This also means that we get a lower bound on $M(n,\Pi_q)$, when $n$ is not the size of any of the candidates.
	
	In the next subsections, we recall the stability results from \cite{BMSzW}. While these results apply to Desarguesian planes, in the introduction, the authors sketch a simple combinatorial proof showing that the minimum number of non-$(r, b)$-colored lines is obtained by coloring $r$ lines entirely red and $b$ lines entirely blue through a common point $P$, whose color may be either red or blue. Assuming that $r\leq b$, the number of non-$(r, b)$-colored lines in such a coloring is either $0$ (if $r=0$ and there were no red points), $1$ (if $r=1$ and only the points of one line are red) or $q+1$ otherwise. This implies that, apart from a few special cases, $M(n,\Pi_q)$ is almost always at least $q+1$.
	
	\subsection{The Desarguesian case for general order}
	Stability results yield exact values and lower bounds for $M^D(n,q)$.
	We do not aim to list all such consequences here, but give a few examples. Here we recall a stability result from \cite{BMSzW} but we reformulate it by using our terminology. 
	
	\begin{result} \cite{BMSzW}
		Let $S$ be a set of points in PG$(2, q)$. Suppose that the number of non-$\kappa(S)$-secants is less than $2q+1$ and they form a dual blocking set. After replacing $S$ by its
		complement if necessary, put
		$r=\kappa(S)\le(q+1)/2$. Then $S$ (or its complement) is one of the following: a single point, a line plus a point, $\kappa(S)$ concurrent lines, $\kappa(S)$ concurrent lines minus their intersection, or a Baer subplane.   
	\end{result}
	
	The difficulty in translating this into a result on $M^D(n, q)$ lies in the specific assumptions involved.
	
	Note that if the non-$\kappa(S)$-secants do not form a dual blocking set, then there exists a point such that all lines through it contain exactly $\kappa(S)$ points of $S$ and hence \[
	|S|\in\bigl\{(q+1)\kappa(S),\,
	1+(q+1)(\kappa(S)-1)\bigr\}.
	\]  The next result  follows from the above result and from Lemma \ref{lem: trivial M} and from Lemma \ref{lem:0,1uj}.
	
	\begin{theorem} \label{cor:des_bounds}
		In the Desarguesian plane PG$(2, q)$,
		the values for $n\in\{0,q+1,q^2, q^2+q+1\}$ are given in Lemma~\ref{lem: n=0, 1}. For $0<n<q^2+q+1$ with $n\notin\{q+1,q^2\}$, the following values and bounds hold for $M^D(n,q)$:
		\begin{itemize}
			\item If $n \equiv 0, 1 \pmod q$, then $M^D(n, q) = q+1$.
			\item If $q \ge 3$, then $M^D(q+2, q) = M^D(q^2 - 1, q) = q+2$. 
			\item If $q$ is a square, then $M^D(q + \sqrt q + 1, q) = M^D(q^2 - \sqrt q, q) = q+ \sqrt{q} + 1$.
			\item If $n\not\equiv 0,1\pmod{q+1}$ and none of the preceding cases applies, then $M^D(n,q)\ge 2q+1$.
		\end{itemize}
	\end{theorem}
	
	\subsection{The Desarguesian case when the order is a prime}
	When the order is a prime, there are stronger results which we recall here
	using our terminology.

	\begin{result} \cite{BMSzW}
		\label{res:prime char}
		Let $p\geq19$ be prime, and let $S$ be a point set in $\PG(2,p)$.
		Suppose that the number of non-$\kappa(S)$-secants is at most
		$\max\{3p+1,4p-22\}$. After replacing $S$ by its complement if
		necessary, put $r=\kappa(S)\leq(p+1)/2$. Then one of the following
		holds:
		\begin{itemize}
			\item[(1)] Apart from a point $P$, the set $S$ consists of the
			points on $r$ lines through $P$. The point $P$ may or may not
			belong to $S$.
			\item[(2)] The set $S$ is obtained from a set in (1) by adding
			or deleting at most two points different from $P$.
			\item[(3)] $r=1$, and $S$ consists of a line together with at
			most three points outside it.
			\item[(4)] $S$ or its complement is projectively equivalent to one the Kiss--Somlai configurations
			described in \cite{BMSzW}, obtained from
			$\{(x,y)\in\mathbb F_p^2:y<x\}$ or its affine complement by
			adding a subset of the three special directions at infinity. Here $x$ and $y$ are identified with their representatives in $\{0,\ldots,p-1\}$, and the three special directions are $0$, $1$ and $\infty$.
			\item[(5)] $r=3$, and $S$ is obtained from the union of the
			three sides of a triangle by possibly deleting some of its
			vertices.
		\end{itemize}
	\end{result}
	
	As before, from this and from Lemma \ref{lem: trivial M} and \ref{lem: S small}, and the constructions above, we get the following corollary for $M^D(n, p)$.
	\begin{corollary}
		Let $p\ge19$ be prime, $N=p^2+p+1$, and put
		$m=\min\{n,N-n\}$. In $\PG(2,p)$, the following
		exceptional exact values hold:
		\[
		\begin{array}{c|l}
			M^D(n,p) & m \\ \hline
			0    & 0\\
			1    & p+1\\
			p+2  & p+2\\
			2p   & 2p+2\\
			2p+2 & p+3\\
			3p-2 & 2p+3,\quad
			\dfrac{p(p-1)}2+3,\quad
			\dfrac{p(p+1)}2-2\\
			3p   & 3p-3\\
			3p+1 & p-2,\quad p+4
		\end{array}
		\]
		For all remaining values of $m$, we have
		\[
		M^D(n,p)=
		\begin{cases}
			p+1,  & m\equiv 0,1\pmod p,\\
			2p+1, & m\equiv -1,2\pmod p,\\
			3p,   & m\equiv -2,3\pmod p,
		\end{cases}
		\]
		while in every other case,
		\[
		M^D(n,p)>\max\{3p+1,4p-22\}.
		\]
	\end{corollary}
	The next result characterizes point sets with small typical intersection numbers. 
	
	\begin{result}
		\label{res:mainprime} \cite{BMSzW}
		Let $S$ be a point set in PG$(2, p)$, $p$ prime. Assume that the number of non-$\kappa(S)$-secants
		of $S$ is less than $d\sqrt p(p+1)$ and suppose also that $\kappa(S) < c\sqrt p$, where $c, d > 0$ and $8d(c+d)  + \frac{20c+28d}{\sqrt p} + \frac{19}{p}< 1$ hold.
		Then we can add and remove in total at most $\lceil d\sqrt{p}\rceil$ points so that we obtain 
		$\kappa(S)$ lines. Furthermore, the size of the union of the set of added or removed points and the intersection points of the $\kappa(S)$ lines is at most $\lceil d\sqrt{p}\rceil$.
	\end{result}
	
	The typical intersection is small, if the size of $S$ is small. Consequently, if $n$ is less than $O(\sqrt p p)$, the stability result in Result \ref{res:mainprime} yields exact values for $M^D(n, p)$ when $n$ matches the size of one of the rigid configurations described above. In particular, it shows the optimality of the construction detailed in Proposition \ref{upperb} (iii).  In all other cases within this range, the result provides a lower bound of $\Omega( p^{\frac32})$ for $M^D(n, p)$. By Lemma \ref{lem: trivial M}, these characterizations and lower bounds naturally extend to sets of large cardinality. 
	
	\begin{result} \cite{BMSzW}
		Consider the set $S$ of points of ${\rm PG}(2, p)$, and assume that the most frequent intersection size of $S$, $\kappa(S)$, is less than $c\sqrt{p}$ while the number of non-$\kappa(S)$-secants of $S$ is less than $d\sqrt p(p+1)$ for some $c, d > 0$ such that the condition $8d(c+d)+\frac{20c+28d}{\sqrt p}+\frac{19}{p}<1$ is satisfied. 
		
		If $\kappa(S) > \lceil d\sqrt p\rceil$, then we can modify (add or delete) in total at most $\lceil d\sqrt p\rceil$ points, so that $S$ becomes the union of $\kappa(S)$ concurrent lines passing through $P$.
	\end{result}
	
	This result and the pencil constructions in Proposition
	\ref{upperb} imply the following corollary for $M^D(n, p)$.
	
	\begin{corollary}
		Let us take the Desarguesian plane $\PG(2, p)$ of prime order $p$.
		Suppose that the integers $r$ and $t$ satisfy $r < c\sqrt{p}$, $1 \le t < d \sqrt p$, where $c,d>0$ and
		$8d(c+d)+\frac{20c+28d}{\sqrt p}+\frac{19}{p}<1$. If $r > \lceil d\sqrt p\rceil$, then
		\[
		M^D(rp+t, p) = tp + 1 - {t-1 \choose 2}.
		\]
	\end{corollary}

	\section{Line sets in general position on level 2}
	\label{sec:sidon}
	
	Construction~4.1(viii) and Lemma~\ref{gen_pos} suggest
	that efficient constructions for $n=qr-t$, with $t>0$
	small, arise from unions of lines after deleting
	suitably chosen intersection points.
	
	In particular, the resulting upper bound on $M(n,\Pi_q)$
	depends on the number of lines meeting the set of pairwise
	intersections. Additional collinear triples among these
	intersection points weaken this bound. This motivates
	the following strengthening of general position.
	
	\begin{defi}[Lines or points in general position on level 2]
		We say that $r$ lines $\ell_1,\ldots,\ell_r$ are in
		\emph{general position on level two} if they are in
		general position and, writing
		\[
		X:=\{\ell_i\cap\ell_j:1\le i<j\le r\},
		\]
		no line other than $\ell_1,\ldots,\ell_r$ contains
		more than two points of $X$.
		
		Dually, we say that $r$ points $P_1,\ldots,P_r$ are
		in \emph{general position on level two} if they are
		in general position and, writing
		\[
		Y:=\{\langle P_i,P_j\rangle:1\le i<j\le r\},
		\]
		no point other than $P_1,\ldots,P_r$ is incident
		with more than two lines of $Y$.
	\end{defi}
	
	\begin{problem}\label{gene_problem}
		Determine
		\[
		\gen_2(q):=
		\max\bigl\{|\mathcal P|:
		\mathcal P\subseteq\PG(2,q)
		\text{ is in general position on level two}\bigr\}.
		\]
	\end{problem}
	
	\begin{theorem}\label{gen_statement}
		We have
		\[
		\gen_2(q)=O(q^{1/2})
		\qquad\text{and}\qquad
		\gen_2(q)=\Omega((q\log q)^{1/5}).
		\]
	\end{theorem}
	
	To obtain lower bounds on $\gen_2(q)$, it is natural
	to restrict the point set to a conic. This leads to
	a projective Sidon-type condition over finite fields.
	
	Consider the conic $X^2=YZ$ in $\PG(2,q)$,
	and parameterize its affine points by
	\[
	P(u)=[u:u^2:1],\qquad u\in\F_q.
	\]
	For distinct $u,v\in\F_q$, the secant $P(u)P(v)$
	has equation
	\[
	(u+v)X-Y-uvZ=0.
	\]
	Consequently, a projective point $Q=[x:y:z]$
	lies on this secant if and only if
	\[
	zuv-x(u+v)+y=0.
	\]
	Moreover, $Q$ lies on the conic precisely when
	$x^2-yz=0$.
	
	For $A\subseteq\F_q$, put
	\[
	H_A:=\{P(u):u\in A\}.
	\]
	Since $H_A$ is contained in a conic, it is
	automatically in general position. It is in general
	position on level two if and only if, for every
	$[a:b:c]\in\PG(2,q)$ with
	\[
	b^2-ac\ne0,
	\]
	the equation
	\[
	auv+b(u+v)+c=0
	\]
	has at most two unordered solutions
	$\{u,v\}\subseteq A$ with $u\ne v$.
	
	This motivates the following arithmetic problem.
	
	\begin{problem}[Projective Sidon sets]
		\label{prob:projective_sidon}
		For an integer $h\ge1$, call $A\subseteq\F_q$
		a \emph{projective $h$-Sidon set} if, for every
		$[a:b:c]\in\PG(2,q)$ satisfying $b^2-ac\ne0$,
		the equation
		\[
		auv+b(u+v)+c=0
		\]
		has at most $h$ unordered solutions
		$\{u,v\}\subseteq A$ with $u\ne v$.
		Determine the maximum size of a projective
		$h$-Sidon set.
	\end{problem}
	
	Thus $H_A$ is in general position on level two
	if and only if $A$ is a projective $2$-Sidon set.
	The condition unifies two familiar types of Sidon
	constraints. Indeed, if $a\ne0$, the equation
	can be rewritten as
	\[
	\left(u+\frac ba\right)
	\left(v+\frac ba\right)
	=\frac{b^2-ac}{a^2},
	\]
	which is a translation-multiplicative condition.
	If $a=0$, then $b\ne0$ and the equation becomes
	\[
	u+v=-\frac cb,
	\]
	an additive Sidon-type condition. The latter case
	corresponds to concurrency at infinity. Hence the
	two conditions are simply coordinate charts of
	a single projective formulation.
	
	We now turn to the proof of Theorem~\ref{gen_statement}.
	For the lower bound we use the following result.
	
	\begin{theorem}[Kostochka--Mubayi--Verstra\"ete
		\cite{KMV}]\label{indep}
		Fix $r\ge2$. There exists a constant $c_r>0$
		such that if $H$ is an $(r+1)$-uniform hypergraph
		on $n$ vertices with maximum $r$-degree
		\[
		0< d<\frac{n}{(\log n)^{3r^2}},
		\]
		then
		\[
		\alpha(H)\ge
		c_r\left(\frac nd\log\frac nd\right)^{1/r}.
		\]
	\end{theorem}
	
	\begin{remark}
		An earlier result of Duke, Lefmann, and
		R\"odl~\cite{DLR} also yields the same order
		of magnitude.
	\end{remark}
	
	\begin{proof}[Proof of Theorem~\ref{gen_statement}]
		For the upper bound, let
		$\mathcal P=\{P_1,\ldots,P_r\}$ be in general
		position on level two. There are
		\[
		\frac12\binom r2\binom{r-2}{2}
		=3\binom r4
		\]
		unordered pairs of secants determined by
		four distinct points of $\mathcal P$.
		
		Their intersection points lie outside
		$\mathcal P$, and the level-two condition
		ensures that these intersection points are
		all distinct. Consequently,
		\[
		3\binom r4\le q^2+q+1-r,
		\]
		which implies $r=O(q^{1/2})$.
		
		For the lower bound, consider the $6$-uniform
		hypergraph $\mathcal H$ with vertex set
		$V=\F_q$. A six-element subset of $V$ is
		an edge if its elements can be partitioned
		into three pairs such that the corresponding
		three secants of the conic are concurrent
		in $\PG(2,q)$.
		
		We claim that
		\[
		\Delta_5(\mathcal H)\le15.
		\]
		Indeed, given five vertices, there are five
		choices for the vertex belonging to the pair
		containing the sixth vertex, and three ways
		to partition the remaining four vertices
		into two pairs. These two pairs determine
		two secants and hence their intersection
		point. The line joining this point to the
		fifth vertex meets the conic in at most one
		additional point. Thus each choice yields
		at most one possible sixth vertex, proving
		the claim.
		
		Applying Theorem~\ref{indep} with $r=5$
		and bounded maximum $5$-degree gives
		\[
		\alpha(\mathcal H)
		=\Omega((q\log q)^{1/5})
		\]
		for sufficiently large $q$.
		
		An independent set $A\subseteq V$ gives
		a point set $H_A$ in general position on
		level two: three concurrent secants at an
		external point would necessarily have
		disjoint pairs of endpoints and hence
		determine an edge of $\mathcal H$.
		Therefore
		\[
		\gen_2(q)\ge\alpha(\mathcal H)
		=\Omega((q\log q)^{1/5}),
		\]
		as required.
	\end{proof}
	
	\begin{remark}
		The hypergraph above satisfies
		$\Delta_k(\mathcal H)=O(q^{5-k})$
		for $1\le k\le5$. These codegree bounds
		are consistent with the scaling of a sparse
		random $6$-uniform hypergraph, but do not
		by themselves yield a matching upper bound
		on its independence number. Determining
		whether the lower bound in
		Theorem~\ref{gen_statement} is sharp
		remains an interesting question.
	\end{remark}

	\noindent \textbf{AI declaration.} The authors   used ChatGPT during drafting and checking of parts of the exposition, and for creating Figure 1.
	
\end{document}